\documentclass{compositio}

\usepackage{amsmath}
\usepackage{amssymb}
\usepackage{mykomma}

\newtheorem*{thm}{Theorem}
\newtheorem*{goddard}{Goddard's Uniqueness Theorem}
\newtheorem*{dong}{Dong's Lemma}
\newtheorem*{ex thm}{Existence Theorem}
\newtheorem{prop}{Proposition}
\newtheorem*{lem}{Lemma}
\newtheorem*{cor}{Corollary}
\theoremstyle{definition}
\newtheorem*{defi}{Definition}
\theoremstyle{remark}
\newtheorem*{conven}{Conventions}

\begin{document}

\title{OPE-Algebras and their Modules}
\author{Markus Rosellen}
\email{rosellen@math.su.se}
\address{Matematiska institutionen\\ Stockholms universitet\\ 
106\,91 Stockholm\\ Sweden}

\thanks{}

\begin{abstract}
Vertex algebras formalize the subalgebra of holomorphic
fields of a conformal field theory. 
OPE-algebras were proposed as a generalization of vertex algebras that 
formalizes the algebra of all fields of a conformal field theory. 
We prove some basic results about them: 
The state-field correspondence is an OPE-algebra isomorphism and 
Dong's lemma and the existence theorem hold for multiply local OPE-algebras; 
locality implies skew-symmetry; 
if skew-symmetry holds then duality implies locality for modules and
they are equivalent for algebras.
We define modules over OPE-algebras. 
\end{abstract}

\maketitle

\section{Introduction}
\label{S:intro}

Conformal field theory (CFT) and vertex algebras made possible dramatic
progress in moonshine, representation theory, quantum topology, moduli spaces,
orbifolds, mirror symmetry, and the geometric Langlands program.
Recently, Kapustin and Orlov \cite{kapustin.orlov.vertex.algebras}
proposed a mathematical definition of CFTs generalizing vertex algebras. 
We call their notion OPE-algebras. ``OPE'' stands for
operator product expansion.

Kapustin and Orlov's paper is about mirror symmetry.
They construct an OPE-algebra $V(T)$ from a (complex) torus $T$ with
a flat metric and a constant $2$-form (a $B$-field).
This is the sigma model of $T$.
The lack of injectivity of the correspondence $T\mapsto V(T)$ is 
called $T$-duality.
In the context of $N=2$ superconformal OPE-algebras,
Kapustin and Orlov determine under which conditions $V(T_1), V(T_2)$ 
are isomorphic and when they are mirror to each other.

OPE-algebras may be defined in terms of a space of fields or
in terms of a state-field correspondence, see section \ref{S:opea}.
The former approach is closer to the Wightman axioms and may apply
to massive deformations of CFTs as well.

According to Polyakov and Kadanoff, the space of fields 
of a CFT endowed with the OPE-coefficients as products is an 
algebra with infinitely many multiplications.
In section \ref{S:alg fields} we define infinitely many products
of two local fields using the OPE.
In section \ref{S:mult loc} we prove Dong's lemma for multiple locality and 
use it to prove the existence theorem and the fact that the state-field
correspondence $Y$ is an OPE-algebra isomorphism.
Thus $Y$ is the adjoint representation.
The existence theorem provides an OPE-algebra structure on a vector space $V$
once a generating set of multiply local fields on $V$ is given.
For vertex algebras, these results are due to Li \cite{li.localsystems},
Lian, Zuckerman \cite{lian.zuckerman.classicalq}, Kac \cite{kac.beginners},
Frenkel, Kac, Radul, Wang \cite{frenkel.kac.radul.wang.w},
and Meurman, Primc \cite{meurman.primc.annihilating.fields.sl2.combinat.ids}.

OPE-algebras are defined in terms of locality.
Borcherds \cite{borcherds.voa} originally defined vertex algebras
in terms of the associativity formula and skew-symmetry.
Li \cite{li.localsystems}
proved the equivalence of these two formulations for vertex algebras.
The associativity formula is the statement that the
state-field correspondence $Y$ is a vertex algebra morphism.
It implies duality. 
Geometrically, the map $Y$ corresponds to 
a 3-punctured sphere and locality and duality are two identities between
the three ways of cutting a 4-punctured sphere into two 3-punctured
spheres \cite{huang.book}. 
Vertex algebra modules can be equivalently defined in terms of the 
associativity formula or duality or the Jacobi identity.

In section \ref{S:skew} we define skew-symmetry for OPE-algebras and
prove that locality implies skew-symmetry.
In section \ref{S:dual} we define duality and prove for algebras
that if skew-symmetry holds then duality is equivalent to locality. 
In section \ref{S:dual skew} we prove for modules that
skew-symmetry and duality imply locality. 
In section \ref{S:mod}
we define the notion of a module over an OPE-algebra $V$ and prove that
if $V$ is uniformly local then $V$ is a $V$-module.
For vertex algebras, these results are due to Li \cite{li.localsystems}.
As in the case of vertex algebras, OPE-algebra modules are interesting 
for example because of their relation to modular invariance.

It is easy to see that any vertex algebra is an OPE-algebra.
Kapustin and Orlov prove the non-trivial result that the
subspace of holomorphic states of an OPE-algebra is a vertex algebra.

Right now, the OPE-algebras $V(T)$ are the only examples
of OPE-algebras besides vertex algebras.
These examples are very special and belong to the class of
{\it additive} OPE-algebras, see \cite{rosellen.addi.opea}.
Additive OPE-algebras are multiply and uniformly local and can be defined
in terms of a ($5$-term) Jacobi identity.
Further examples of OPE-algebras should exist, e.g.~or\-bifolds of $V(T)$
and various rational models (minimal models, WZW-models etc.).
For their construction one should use results about intertwiners,
e.g.~the fundamental results of Huang,
see \cite{huang.differ.equations.intertwining.operators}.

\begin{conven}
We work over a field $\K$ of characteristic $0$.
We always denote by $E$ a vector space and by
$A$ an associative algebra.
For $a\in A$ and $n\in\K$, define
$a^{(n)}:=a^n/n!$ if $n\in\N:=\Z_{\geq 0}$ and $a^{(n)}:=0$ otherwise.
Let $[n]:=\set{1, \dots, n}$ for $n\in\N$.

We work with super objects without making
this explicit in our terminology,
e.g. super vector spaces are just called vector spaces.
Linear maps need not be even, i.e.~preserve the super grading.
Supersigns are written as powers of $\zeta:=-1$.
The parity of $a$ is denoted by $\ta\in\set{\bar{0},\bar{1}}$.
\end{conven}

\section{OPE-Algebras}
\label{S:opea}

We define OPE-algebras in terms of a space of fields and show that
this definition is equivalent to a definition in terms of a
state-field correspondence.

\medskip

An $E$-valued {\it distribution} 
is a formal sum $a(z)=\sum_{n\in\K}a_n z^{-n-1}$ where 
$z$ is a formal variable and $a_n\in E$. 
Here and in the following we only discuss the case of one
variable. 

The vector space $E\set{z}$ of distributions
is a module over the group ring $\K[z^{\K}]$ of $\K$,
$z^m a(z):=\sum a_{n+m}z^{-n-1}$. 
Let $E\sqbrack{z}$ be the submodule generated
by the subspace of power series 
$E\pau{z}:=\set{\sum_{n\in\Z_{<0}}a_n z^{-n-1}}$. 
There exists a morphism
$E\sqbrack{z,w}\to E\sqbrack{z}, a(z,w)\mapsto a(z,z)$.
A linear map $E\otimes F\to G$ induces morphisms 
$E\set{z}\otimes F\set{w}\to G\set{z,w}$ and
$E\sqbrack{z}\otimes F\sqbrack{z}\to G\sqbrack{z}$.
Define $\del_z a(z):=-\sum n\, a_{n-1}z^{-n-1}$.

Let $z, \bz$ be variables. 
We shall use the following notations time and again.
Define $a(\vz):=a(z,\bz)$ and $E\set{\vz}:=E\set{z,\bz}$.
If $S, T$ are sets and $\vs\in S\times T$ then 
$s, \bs$ denote the first and the second component of $\vs$.
If $a_s(z), \ba_{\bs}(z)\in A\set{z}$ 
then we define $\va_{\vs}(\vz):=a_s(z)\ba_{\bs}(\bz)$.
For example,
if $\vn\in\K^2$ and $\va\in A^2$ then
$\vz^{\vn}=z^n \bz^{\bn}$ and $\va^{(\vn)}=a^{(n)}\ba^{(\bn)}$.
Thus $a(\vz)=\sum_{\vn\in\K^2}a_{\vn}\vz^{-\vn-1}$
with $1:=(1,1)\in\K^2$. 

Unless stated otherwise, from now on all distributions will be
$\End(E)$-valued.

A distribution $a(\vz)$ is a {\it field} if
$a(\vz)b\in E\sqbrack{\vz}$ for any $b\in E$. 
Let $\cF_r(E)$ be the space of fields 
$a(\vz_1, \dots, \vz_r)$ and $\cF(E):=\cF_1(E)$.

Let $1\in E\even$ and $\vT\in\End(E)\even^2$.
We call $1$ {\it invariant} if $T1=\bT 1=0$.
Define $s_1:\End(E)\set{\vz}\to E$ by $a(\vz)\mapsto a_{-1}(1)$. 
A distribution $a(\vz)$ is {\it weakly creative} for $1$ if
$a(\vz)1\in E\pau{\vz}$.
It is {\it creative} for $1$ and $\vT$ if
$a(\vz)1=e^{\vz\vT}s_1 a(\vz)$.
It is {\it translation covariant} for $\vT$
if $[T,a(\vz)]=\del_z a(\vz)$ and $[\bT,a(\vz)]=\del_{\bz} a(\vz)$.
A subspace $\cF\subset\End(E)\set{\vz}$ is {\it complete}
if $s_1 \cF=E$.

If $1$ is invariant and $a(\vz)$ translation covariant and weakly creative
then $a(\vz)$ is creative.
If $[T,\bT]=0$ then translation covariance is equivalent to
$e^{\vw\vT}a(\vz)e^{-\vw\vT}=a(\vz+\vw)$ where 
$a(z+w):=e^{w\del_z}a(z)$.

For $\vn\in\K^2$, define
$$
(\vz-\vw)^{\vn}
\; :=\;
\sum_{\vi\in\N^2}(-1)^{\vi}\binom{\vn}{\vi}\, \vz^{\vn-\vi}\vw^{\vi}.
$$
For $\vn\in\vbbK:=\set{\vn\in\K^2\mid n-\bn\in\Z}$, define
$(-1)^{\vn}:=(-1)^{n-\bn}$ and 
$(\vz-\vw)_{w>z}^{\vn}:=(-1)^{\vn}(\vw-\vz)^{\vn}$.

Distributions $a(\vz), b(\vz)$ are {\it local}
if there exist $c^i(\vz,\vw)\in\cF_2(E)$ and $\vh_i\in\vbbK$
such that
\begin{equation}
\label{E:ope}
a(\vz)b(\vw)
\; =\;
\sum_{i=1}^r \frac{c^i(\vz,\vw)}{(\vz-\vw)^{\vh_i}},
\qquad
\paraab\, b(\vw)a(\vz)
\; =\;
\sum_{i=1}^r \frac{c^i(\vz,\vw)}{(\vz-\vw)_{w>z}^{\vh_i}}.
\end{equation}
Equations \eqref{E:ope} are called {\it OPEs} in $z>w$ and $w>z$.

If $\cP$ is a property of elements or pairs of elements of a
set $S$ then we say that a subset $T\subset S$ satisfies $\cP$
if $\cP$ is satisfied for any element, resp., any pair of elements of $T$.

\begin{defi}
A vector space $V$ together with a vector $1\in V\even$ and a 
subspace $\cF\subset\End(V)\set{\vz}$ is an {\it OPE-algebra} 
if there exists $\vT\in\End(V)\even^2$ such that 
$1$ is invariant and 
$\cF$ is weakly creative, translation covariant, complete, and local.
\end{defi}

The following result is proven in \cite{kapustin.orlov.vertex.algebras}.

\begin{goddard}
Let $\cF\subset\End(E)\set{\vz}$ be a creative, complete, local subspace.
Then $s_1:\cF\to E$ is an isomorphism. \hfill $\square$
\end{goddard}

Let $V$ be an OPE-algebra. The theorem shows that 
the inverse $Y:=(s_1|_{\cF})\inv: V\to\End(V)\set{\vz}$ exists
and if $a(\vz)\in\End(V)\set{\vz}$ is creative and local to $\cF$ then 
$a(\vz)=Y(s_1 a(\vz))\in\cF$. 

Let $V$ be a vector space.
To give an even linear map $Y:V\to\End(V)\set{\vz},
a\mapsto a(\vz)=\sum a_{(\vn)}\vz^{-\vn-1}$, is equivalent to
giving an even multiplication
$V\otimes V\to V, a\otimes b\mapsto a_{(\vn)}b$, for any $\vn\in\K^2$.
We call this a {\it $\K^2$-fold algebra}. 

Let $V$ be a $\K^2$-fold algebra, $1\in V\even$, and $\vT\in\End(V)\even^2$.
Define $\cF_Y:=Y(V)$.
Define $\vT_1\in\End(V)\even^2$ by
$T_1a:=a_{(-2,-1)}1$ and $\bT_1a:=a_{(-1,-2)}1$.
We call $1$ a {\it (weak) right identity}
if $\cF_Y$ is (weakly) creative for $1$ and $\vT_1$ and
$s_1 a(\vz)=a$ for any $a\in V$. 

We call $1$ a {\it left identity}
if $Y(1)=1(z)$ where $1(z):=\id_V$ is the {\it identity field}.
An {\it identity} is a left and right identity.
Assume that $[T,\bT]=0$.
The pair $\vT$ is a {\it translation generator}
if $\cF_Y$ is translation covariant for $\vT$.
It is a {\it translation endomorphism}
if $(Ta)(\vz)=\del_z a(\vz)$ and $(\bT a)(\vz)=\del_{\bz}a(\vz)$.

To give a vector space $V$, a vector $1\in V\even$,
and a weakly creative subspace $\cF\subset\End(V)\set{\vz}$ 
such that $s_1:\cF\to V$ is an isomorphism 
is equivalent to giving a $\K^2$-fold algebra with 
a weak right identity.
Moreover, 
Proposition \ref{P:skew} shows that OPE-algebras are unital
and thus a weak right identity is unique.
This shows that the above definition is equivalent to the following one.

\begin{defi}
An {\it OPE-algebra} is a $\K^2$-fold algebra such that there exist
a translation generator $\vT$ and an invariant weak right identity $1$
and $\cF_Y$ is local.
\end{defi}

It follows that $\vT=\vT_1$ and that $\vT$ is a translation endomorphism.
Moreover, one need not require that $[T,\bT]=0$, it is a consequence.
Proposition \ref{P:skew} shows that $1$ is an identity.
A {\it morphism} of OPE-algebras is a morphism of the underlying
unital $\K^2$-fold algebras.

If $V, W$ are OPE-algebras then so is $V\otimes_{\K} W$.

\section{The Algebra of Fields}
\label{S:alg fields}

We define a field $a(\vz)_{(\vn)}b(\vz)$ for any local distributions
$a(\vz), b(\vz)$ and any $\vn\in\K^2$ and prove some basic properties of 
this (partial) $\K^2$-fold algebra.

\medskip

Distributions in $E\pau{z\uppm}:=\{\sum_{n\in\Z}a_n z^{-n-1}\}$
are called {\it holomorphic}.
Define $\cF_z(E):=\cF(E)\cap\End(E)\pau{z\uppm}$. 
Define
$(z-w)^n [a(z),b(\vw)]:=(z-w)^n a(z)b(\vw)-\paraab(z-w)_{w>z}^n b(\vw)a(z)$
for $a(z)\in\cF_z(E), b(\vz)\in\cF(E)$, and $n\in\Z$.
The following two results are proven in \cite{kapustin.orlov.vertex.algebras}.

\begin{prop}
\label{P:ko}
\begin{enumerate}
\item
Let $c^i(\vz,\vw)\in E\sqbrack{\vz,\vw}$ and $\vh_i\in\K^2$ 
such that $\vh_i\notin \vh_j+\Z^2$ for $i\ne j$ and
$$
\sum_{i=1}^r\:
c^i(\vz,\vw)\, (\vz-\vw)^{\vh_i}
\; =\;
0.
$$
Then $c^i(\vz,\vw)=0$ for any $i$.

\item
Let $a(z)\in\cF_z(E)$ and $b(\vz)\in\cF(E)$. 
Then $a(z), b(\vz)$ are local iff 
there exists $N\in\Z$ such that $(z-w)^N [a(z),b(\vw)]=0$.
In this case we have the OPEs
$$
a(z)b(\vw)
\; =\;
\frac{c(z,\vw)}{(z-w)^N},
\qquad
\paraab b(\vw)a(z)
\; =\;
\frac{c(z,\vw)}{(z-w)_{w>z}^N}.
$$\hfill $\square$
\end{enumerate}
\end{prop}

An OPE \eqref{E:ope} is {\it reduced} if 
$\vh_i\notin \vh_j+\Z^2$ for $i\ne j$.
Reduced OPEs always exist for local distributions.

Let $a(\vz), b(\vz)\in\End(E)\set{\vz}$ be local with OPE
\eqref{E:ope}. For $\vn\in\K^2$, define
$$
a(\vw)_{(\vn)}b(\vw)
\; :=\;
\sum_{i=1}^r\:
\del_{\vz}^{(\vh_i-1-\vn)}c^i(\vz,\vw)|_{\vz=\vw}.
$$
Proposition \ref{P:ko}\,i) and $\del_{\vz}^{(\vn+\vi)}
((\vz-\vw)^{\vi}c(\vz,\vw))|_{\vz=\vw}
=\del_{\vz}^{(\vn)}c(\vz,\vw)|_{\vz=\vw}, \vi\in\N^2$,
show that this definition does not depend on
the choice of the OPE.
The fields $a(\vz)_{(\vn)}b(\vz)$ are the Taylor coefficients
of a non-existing Taylor expansion of $c^i(\vz,\vw)$.
If \eqref{E:ope} is reduced then
$c^i(\vz,\vz)=a(\vz)_{(\vh_i-1)}b(\vz)$.

If $a(z)\in\cF_z(E)$ and $b(\vz)\in\cF(E)$ then 
$$
a(w)_{(\vn)}b(\vw)
\; =\;
a(w)_{(n)}b(\vw)
\; :=\;
\res_z(z-w)^n [a(z),b(\vw)]
$$
for $\vn\in\Z\times\set{-1}$ and $a(z)_{(\vn)}b(\vz)=0$ otherwise,
see \cite{kapustin.orlov.vertex.algebras}.
Thus the products $a(\vz)_{(\vn)}b(\vz)$ generalize the
products $a(z)_{(n)}b(\vz)$ known from the theory of vertex algebras
\cite{li.localsystems,lian.zuckerman.classicalq,kac.beginners}.
Note that $a(z)_{(n)}b(\vz)$ is always defined
whereas $a(\vz)_{(\vn)}b(\vz)$ is only defined if an OPE in $z>w$ exists.

Let $V$ be an OPE-algebra.
It is shown in \cite{kapustin.orlov.vertex.algebras}
that $V_z:=\set{a\in V\mid Ya\in\End(V)\pau{z\uppm}}$
is a {\it vertex subalgebra}, 
i.e.~$V_z\subset V$ is a unital $\K^2$-fold subalgebra,
$a_{(\vn)}b=0$ for $a, b\in V_z, \vn\notin\Z\times\set{-1}$,
and $(V,(_{(n,-1)})_{n\in\Z})$ is a vertex algebra.
This subalgebra is the {\it chiral algebra} of $V$.
Moreover, they show that
$(a_{(n,-1)}b)(\vz)=a(z)_{(n)}b(\vz)$ for any $a\in V_z, b\in V$.
Similarly, the {\it anti-chiral} algebra $V\subbz$ is defined and
$V_z, V\subbz$ commute, i.e.~$[a(z),b(\bz)]=0$.

For $a(z)\in E\set{z}$ and $S\subset\K$, define 
$\supp_z a(z):=\set{n\mid a_n\ne 0}$ and
$a(z)|_S:=\sum_{n\in S}a_n z^{-n-1}$.
Let $E\set{\vz}':=
\set{a(\vz)\mid\supp_{\vz} a(\vz)\subset\vbbK}$.
A $\K^2$-fold algebra $V$ is a {\it bounded $\vbbK$-fold} algebra
if $\cF_Y\subset\cF(V)\cap\End(V)\set{\vz}'$.

\begin{prop}
\label{P:algebra}
Let $a(\vz), b(\vz)\in\End(E)\set{\vz}$ be local. 

\begin{enumerate}
\item
$(\del_z,\del_{\bz})$ is a translation endomorphism and 
a translation generator for $a(\vz), b(\vz)$
(i.e.~we have 
$\del_z a(\vz)_{(\vn)}b(\vz)=-n\, a(\vz)_{(n-1,\bn)}b(\vz)$ etc.).
The identity field  $1(z)$ is an identity for any $c(\vz)\in\cF(E)$
(i.e.~we have $1(z)_{(\vn)}c(\vz)=\de_{\vn,-1}c(\vz)$ etc.).
Moreover, $\vT_{1(z)}=(\del_z,\del_{\bz})$.

\item
If $a(\vz), b(\vz)$ are creative and translation covariant 
then so is $a(\vz)_{(\vn)}b(\vz)$ and 
$s_1(a(\vz)_{(\vn)}b(\vz))=a_{\vn}s_1b(\vz)$ for any $\vn$. 

\item
OPE-algebras are bounded $\vbbK$-fold algebras and 
$c^i(\vz,\vw)\in\End(E)\set{\vz,\vw}'$.
\end{enumerate}
\end{prop}

\begin{proof}
Part i) and the fact that $a(\vz)_{(\vn)}b(\vz)$ is translation covariant
follow from a direct calculation.

Let \eqref{E:ope} be a reduced OPE of $a(\vz), b(\vz)$ and $c\in E$.
Define $S:=\supp_{\vz} b(\vz)c+\Z^2$ and $T:=\K^2\times(\K^2\setminus S)$.
Because $(\vz-\vw)^{\vh}\in\K\set{\vz}\pau{\vw}$ we have
$$
\sum_i\:
\frac{c^i(\vz,\vw)c|_T}{(\vz-\vw)^{\vh_i}}
\; =\;
\sum_i\:
\frac{c^i(\vz,\vw)c}{(\vz-\vw)^{\vh_i}}
\bigg|_T
\; =\;
a(\vz)b(\vw)c|_T
\; =\;
0.
$$
Proposition \ref{P:ko}\,i) yields $c^i(\vz,\vw)c|_T=0$. 
Hence $\supp_{\vw}c^i(\vz,\vw)c\subset S$.
Similarly, the second OPE yields 
$\supp_{\vz}c^i(\vz,\vw)c\subset \supp_{\vz} a(\vz)c+\Z^2$.
In particular,  $c^i(\vz,\vw)1\in E\pau{\vz\uppm,\vw\uppm}$.

Define $S_i:=(\vh_i+\Z^2)\times\K^2$. 
Because $(\vz-\vw)^{\vh}\in\vz^{\vh}\K\pau{\vz\uppm}\set{\vw}$ we have
$$
a(\vz)b(\vw)1|_{S_j}
=
\sum_i
\frac{c^i(\vz,\vw)1}{(\vz-\vw)^{\vh_i}}
\bigg|_{S_j}
=
\frac{c^j(\vz,\vw)1}
{(\vz-\vw)^{\vh_j}}.
$$
Since $b(\vw)1\in E\pau{\vw}$ we get $c^j(\vz,\vw)1\in E\set{\vz}\pau{\vw}$.
Similarly, the second OPE yields $c^i(\vz,\vw)1\in E\set{\vw}\pau{\vz}$.
Thus $c^i(\vz,\vw)1\in E\pau{\vz,\vw}$.
This implies that $a(\vz)_{(\vn)}b(\vz)$ is weakly creative.
Applying the reduced OPE to $1$ and setting $\vw=0$ yields
$$
a(\vz)s_1b(\vz)
\; = \;
\sum_i\:
\frac{c^i(\vz,\vw)1|_{\vw=0}}{\vz^{\vh_i}}
\; \in \;
\sum_i\:
\vz^{-\vh_i}E\pau{\vz}.
$$
In particular, we obtain iii).
Moreover, this shows that for any $\vn\in\N^2$ we have 
$$
a(\vw)_{(\vh_i-1-\vn)}b(\vw)1|_{\vw=0}
\; =\;
\del_{\vz}^{(\vn)} c^i(\vz,\vw)1|_{\vz=\vw=0}
\; =\;
a_{\vh_i-1-\vn}s_1b(\vz).
$$
Hence $s_1(a(\vz)_{(\vn)}b(\vz))=a_{\vn}s_1b(\vz)$ for any $\vn\in\K^2$.
Note that $c(\vz)$ is creative 
iff $c(\vz)$ is weakly creative and 
$s_1\del_{\vz}^{(\vn)}c(\vz)=\vT^{(\vn)}s_1 c(\vz)$ for any $\vn$. 
By i) we know that 
$\del_{\vz}$ is a translation generator for local distributions.
Thus 
\begin{align}
\notag
s_1\del_{\vz}^{(\vm)}(a(\vz)_{(\vn)}b(\vz))
&=
s_1\sum_{\vi=0}^{\vm}
(-1)^{\vi}\binom{\vn}{\vi} a(\vz)_{(\vn-\vi)}b(\vz)
+
a(\vz)_{(\vn)}\del_{\vz}^{(\vm-\vi)}b(\vz)
\\
\notag
&=
\sum_{\vi=0}^{\vm}
(-1)^{\vi}\binom{\vn}{\vi} a_{\vn-\vi}s_1b(\vz)
+
a_{\vn}s_1\del_{\vz}^{(\vm-\vi)}b(\vz).
\end{align}
On the other hand,
since $a(\vz)$ is translation covariant we have 
\begin{align}
\notag
\vT^{(\vm)}s_1(a(\vz)_{(\vn)}b(\vz))
\; &=\;
\vT^{(\vm)}a_{\vn}s_1b(\vz)
\\
\notag
&=\;
\sum_{\vi=0}^{\vm}\:
(-1)^{\vi}\binom{\vn}{\vi}\, a_{\vn-\vi}s_1b(\vz)
\, +\,
a_{\vn}\vT^{(\vm-\vi)}s_1b(\vz).
\end{align}
This shows that $a(\vz)_{(\vn)}b(\vz)$ is creative.
\end{proof}

\section{Multiple Locality}
\label{S:mult loc}

We prove Dong's lemma, the existence theorem, and 
that $Y$ is an OPE-algebra isomorphism for multiply local fields.

\medskip

For a permutation $\si\in\bbS_r$, $\vn\in\vbbK$, and $i\ne j\in [r]$, define
$$
(\vz_i-\vz_j)_{\si}^{\vn}
\; :=\;
\begin{cases}
\qquad (\vz_i-\vz_j)^{\vn}&
\text{if $\si^{-1}(i)<\si^{-1}(j)$ and}\\
\qquad (\vz_i-\vz_j)_{z_j>z_i}^{\vn}&
\text{otherwise.}
\end{cases}
$$

Distributions $a^1(\vz), \dots, a^r(\vz)\in\End(E)\set{\vz}$ 
are {\it multiply local} 
if there exist $r_{ij}\in\N, \vh_{ij}^k\in\vbbK$, and $c^{\al}\in\cF_r(E)$ 
for any $i,j\in [r],\: k\in [r_{ij}]$, and $\al\in\prod_{i<j}[r_{ij}]$
such that for any $\si\in\bbS_r$ we have
\begin{equation}
\label{E:mult loc}
\zeta'\,
a^{\si 1}(\vz_{\si 1})\dots  a^{\si r}(\vz_{\si r})
\; =\;
\sum_{\al}\:
\frac{c^{\al}(\vz_1,\dots,\vz_r)}
{\prod_{i<j}(\vz_i-\vz_j)_{\si}^{\vh_{ij}^{\al(i,j)}}}
\end{equation}
where $\zeta'$ is the obvious supersign.
Equations \eqref{E:mult loc} are also called {\it OPEs}.
They are {\it reduced} if 
$\vh_{ij}^k\notin\vh_{ij}^l+\Z^2$ for $k\ne l$.

A pair of distributions is multiply local iff it is local.
A subset $S\subset\End(E)\set{\vz}$ is {\it multiply local} if 
any finite family in $S$ is multiply local.
An OPE-algebra $V$ is {\it multiply local} if $\cF_Y$ is.

A subset $S\subset\End(E)\set{\vz}$ is {\it uniformly local} if 
there exist $r_{ab}\in\N$ and $\vh_{ab}^k\in\vbbK$ 
for any $a(\vz), b(\vz)\in S, k\in [r_{ab}]$ such that
\eqref{E:mult loc} is satisfied for any $a^i(\vz)\in S$ 
with $\vh_{ij}^k:=\vh_{a^i a^j}^k$ and 
$c\upal 1\in E\pau{\vz_1, \dots, \vz_r}$.
For example, additive OPE-algebras are uniformly local,
see \cite{rosellen.addi.opea}.
The proof of Proposition \ref{P:algebra}
shows that $c\upal 1\in E\pau{\vz_1,\vz_2}$ if 
$r=2$ and $a^1(\vz), a^2(\vz)$ are local and weakly creative.

\begin{dong}
\begin{enumerate}
\item
Let $a^1(\vz), \dots, a^r(\vz)$ be multiply local and
$a^1(\vz), a^2(\vz)$ local.
Then $a^1(\vz)_{(\vn)}\linebreak[0]a^2(\vz), a^3(\vz), 
\dots, a^r(\vz)$
are multiply local for any $\vn$.

\item
Let $S\subset\cF(E)$ be multiply (uniformly) local. 
Then there exists a multiply (uniformly) local subspace $\sqbrack{S}$
that is closed
with respect to the products $a(\vz)_{(\vn)}b(\vz)$, contains $1(z)$,
and is generated by $S$ as a unital $\K^2$-fold algebra.
\end{enumerate}
\end{dong}

\begin{proof}
i)\:
We may assume that \eqref{E:mult loc} and 
$a^1(\vz)a^2(\vw)=
\sum_{k=1}^{r_{12}}(\vz-\vw)^{-\vh_{12}^k}\, c^k(\vz,\vw)$
are reduced OPEs.
Let $\si\in\bbS_r$ such that $\si^{-1}(2)=\si^{-1}(1)+1$.
Then
$$
\zeta'\,
a^{\si 1}(\vz_{\si 1})\dots  a^{\si r}(\vz_{\si r})
\; =\;
\sum_{k=1}^{r_{12}}
\frac{a^{\si 1}(\vz_{\si 1})\dots c^k(\vz_1,\vz_2)\dots a^{\si r}(\vz_{\si r})}
{(\vz_1-\vz_2)^{\vh_{12}^k}}.
$$
Multiple locality and Proposition \ref{P:ko}\,i) yield
$$
\zeta'\,
a^{\si 1}(\vz_{\si 1})\dots c^k(\vz_1,\vz_2)\dots a^{\si r}(\vz_{\si r})
\; =\;
\sum_{\al:\: \al(1,2)=k}\:
\frac{c^{\al}(\vz_1,\dots,\vz_r)}
{\prod_{(i,j)\ne (1,2)}(\vz_i-\vz_j)_{\si}^{\vh_{ij}^{\al(i,j)}}}
$$
for any $k$.
We have
$a^1(\vw)_{(\vh_{12}^k-1-\vn)}a^2(\vw)=
\del_{\vz}^{(\vn)}c^k(\vz,\vw)|_{\vz=\vw}$
for any $\vn\in\N^2$.
Thus by acting with $\del_{\vz_1}^{(\vn)}$ on
the last equation and setting $\vz_1=\vz_2$ we obtain
$$
\zeta'\,
a^{\si 1}(\vz_{\si 1})\dots 
a^1(\vz_2)_{(\vh_{12}^k-1-\vn)}a^2(\vz_2)\dots a^{\si r}(\vz_{\si r})
\; =\;
\sum_{\al:\: \al(1,2)=k}\:
\frac{\tc^{\al}(\vz_2,\dots,\vz_r)}
{\prod_{1<i<j}(\vz_i-\vz_j)_{\si}
^{\vh_{ij}^{\al(i,j)}+\de_{i,2}(\vh_{1j}^{\al(1,j)}+\vn)}}
$$
for some fields $\tc^{\al}$.

ii)\:
This follows from i) and its proof.
\end{proof}

\begin{thm}
Let $V$ be a multiply local OPE-algebra. 
Then $\cF_Y$ is closed with respect to the products $a(\vz)_{(\vn)}b(\vz)$ 
and $Y:V\to\cF_Y$ is an OPE-algebra isomorphism.
\end{thm}

\begin{proof}
Dong's lemma and Proposition \ref{P:algebra}\,ii) show that 
$\sqbrack{\cF_Y}$ is local and creative.
Goddard's uniqueness theorem and Proposition \ref{P:algebra}\,ii) imply that
$a(\vz)_{(\vn)}b(\vz)=Y(s_1(a(\vz)_{(\vn)}b(\vz)))=(a_{(\vn)}b)(\vz)$.
We have $Y1=1(z)$ since $V$ is unital by Proposition \ref{P:skew}.
\end{proof}

\begin{ex thm}
Let $V$ be a vector space, $1\in V\even, \vT\in\End(V)\even^2$,
and $S\subset\cF(V)$ a 
weakly creative, translation covariant, multiply (uniformly) local subset 
such that $1$ is invariant and
$$
V
\; =\;
\rspan\set{ a^1_{\vn_1}\dots a^r_{\vn_r}1\mid
a^i(\vz)\in S,\, \vn_i\in\K^2,\, r\in\N }.
$$
Then there exists a unique multiply (uniformly) local
OPE-algebra structure $Y$
on $V$ such that $1$ is a weak right identity and
$Y(s_1 a(\vz))=a(\vz)$ for any $a(\vz)\in S$.
We have $\cF_Y=\sqbrack{S}$.
\end{ex thm}

\begin{proof}
Dong's lemma and Proposition \ref{P:algebra}\,ii) show that 
$\sqbrack{S}$ is creative, translation covariant, complete, 
and multiply (uniformly) local.
Thus $(V,1,\sqbrack{S})$ is an OPE-algebra satisfying
the three properties of the theorem.
If $Y'$ is such an OPE-algebra structure 
then $\sqbrack{S}=\cF_{Y'}$ and $Y'=(s_1|_{\sqbrack{S}})\inv$.
\end{proof}

\section{Skew-Symmetry}
\label{S:skew}

We define skew-symmetry and prove that locality implies skew-symmetry.

\medskip

Let $V$ be a bounded $\vbbK$-fold algebra and $\vT\in\End(V)\even^2$.
{\it Skew-symmetry} for $\vT$ is the identity 
$$
\paraab\, b_{(\vn)}a
\; =\; 
\sum_{\vi\in\N^2}\:
(-1)^{\vn+\vi}\: \vT^{(\vi)}(a_{(\vn+\vi)}b)
$$
for any $\vn\in\vbbK$.
This is equivalent to $\paraab\, b(\vz)a=e^{\vz\vT}a(-\vz)b$
where $c(-\vz):=\sum_{\vn\in\vbbK}(-1)^{\vn}c_{\vn}\vz^{-\vn-1}$
for $c(\vz)\in E\set{\vz}'$.

\begin{prop}
\label{P:skew}
\begin{enumerate}
\item
Local distributions in $\End(E)\set{\vz}$
satisfy skew-symmetry for $(\del_z,\del_{\bz})$.

\item
OPE-algebras satisfy skew-symmetry and are unital.
\end{enumerate}
\end{prop}

\begin{proof}
i)\: 
Let $a(\vz), b(\vz)$ be local with reduced OPE \eqref{E:ope}
and $\vn\in\N^2$. 
The OPE in $w>z$, the binomial formula, and the chain rule yield
\begin{align}
\notag
\paraab\,b(\vw)_{(\vh_i-1-\vn)}a(\vw)
\;=\;
&
(-1)^{\vh_i}\,\del_{\vz}^{(\vn)}c^i(\vw,\vz)|_{\vz=\vw}
\\
\notag
=\;
&
\sum_{\vk\in\N^2}\:
(-1)^{\vh_i}\,
(\del_{\vz}+\del_{\vw})^{(\vk)}(-\del_{\vw})^{(\vn-\vk)}
c^i(\vw,\vz)|_{\vz=\vw}
\\
\notag
=\;
&
\sum_{\vk\in\N^2}\:
(-1)^{\vh_i+\vn+\vk}\: 
\del_{\vw}^{(\vk)}(a(\vw)_{(\vh_i-1-\vn+\vk)}b(\vw)).
\end{align}

ii)\:
Because of Proposition \ref{P:algebra}\,ii), part i), and
creativity of $a(\vz)_{(\vn+\vi)}b(\vz)$ we have
\begin{align}
\notag
\paraab\, b_{(\vn)}a
\; =\;
\paraab\, s_1(b(\vz)_{(\vn)}a(\vz))
\; &=\;
s_1\sum_{\vi\in\N^2}\:
(-1)^{\vn+\vi}\, \del_{\vz}^{(\vi)}(a(\vz)_{(\vn+\vi)}b(\vz))
\\
\notag
&=\;
\sum_{\vi\in\N^2}\:
(-1)^{\vn+\vi}\, \vT^{(\vi)}(a_{(\vn+\vi)}b).
\end{align}
Skew-symmetry implies that $1$ is a left identity iff $1$ is a right
identity.
\end{proof}

\section{Duality and Locality}
\label{S:dual}

We define duality and prove for $\vbbK$-fold algebras that 
if skew-symmetry holds then duality is equivalent to locality.

\medskip 

Let $V$ be a bounded $\vbbK$-fold algebra and $M$ a
bounded $\vbbK$-fold $V$-{\it module},
i.e.~$M$ is a vector space with an even linear map 
$Y: V\to\cF(M)\cap\End(M)\set{\vz}'$.

Elements $a\in V, c\in M$ are {\it dual} in the {\it direct channel}
if there exist $d^j(\vw,\vx)\in\cF_2(V)$ and $\vh_{ac}^j\in\vbbK$
such that 
\begin{equation}
\label{E:dual}
a(\vx+\vw)b(\vw)c
\; =\;
\sum_{j=1}^s\frac{d^j(\vw,\vx)b}{(\vx+\vw)^{\vh_{ac}^j}},
\qquad
(a(\vx)b)(\vw)c
\; =\;
\sum_{j=1}^s\frac{d^j(\vw,\vx)b}{(\vw+\vx)^{\vh_{ac}^j}}
\end{equation}
for any $b\in V$.
Elements $a, b\in V$ are {\it local} on $M$ in the {\it direct channel}
if $a(\vz), b(\vz)\in\cF(M)$ are local.

The term ``direct channel'' means that the singularity in
\eqref{E:ope} and \eqref{E:dual} is the one corresponding to $a, b$, resp.,
$a, c$.

\begin{prop}
\label{P:dual loc}
Let $V$ be a bounded $\vbbK$-fold algebra satisfying skew-symmetry for a
translation generator $\vT$.
Then $a, b\in V$ are local in the direct channel iff
$a, b$ are dual in the direct channel.
\end{prop}

The idea of the proof is to argue that $a(bc)=b(ac)$ iff $a(cb)=(ac)b$.

\begin{proof}
`$\Rightarrow$'\: Let \eqref{E:ope} be an OPE of $a(\vz), b(\vz)$. 
For $c\in V$, we have
$$
\zeta^{\tb\tc}\, a(\vx+\vw)c(\vw)b
\; =\;
a(\vx+\vw)e^{\vw\vT}b(-\vw)c
\; =\;
e^{\vw\vT}a(\vx)b(-\vw)c
\; =\;
\sum_i\frac{e^{\vw\vT}c^i(\vx,-\vw)c}{(\vx+\vw)^{\vh_i}}
$$
and
$$
\zeta^{\tb\tc}\, (a(\vx)c)(\vw)b
\; =\; 
\zeta^{\ta\tb}\, e^{\vw\vT}b(-\vw)a(\vx)c
\; =\; 
\sum_i\frac{e^{\vw\vT}c^i(\vx,-\vw)c}{(\vw+\vx)^{\vh_i}}.
$$

`$\Leftarrow$'\: This is proven in the same way.
\end{proof}

\section{Duality and Skew-Symmetry}
\label{S:dual skew}

We prove for modules that duality and skew-symmetry imply locality.

\medskip 

Let $V$ be a bounded $\vbbK$-fold algebra and 
$M$ a bounded $\vbbK$-fold $V$-module.
Elements $a, b\in V$ are {\it dual} in the {\it exchange channel}
if there exist $\vh_{ab}^i\in\vbbK$ and for any $c\in M$ there exist 
$d^j(\vw,\vx)\in\sum_{i=1}^r \vx^{-\vh_{ab}^i} M\pau{\vx}\sqbrack{\vw}$ 
and $\vh_{ac}^j\in\vbbK$ such that 
\begin{equation}
\label{E:dual ab}
a(\vx+\vw)b(\vw)c
\; =\;
\sum_{j=1}^s\frac{d^j(\vw,\vx)}{(\vx+\vw)^{\vh_{ac}^j}},
\qquad
(a(\vx)b)(\vw)c
\; =\;
\sum_{j=1}^s\frac{d^j(\vw,\vx)}{(\vw+\vx)^{\vh_{ac}^j}}.
\end{equation}

\begin{prop}
\label{P:dual skew}
Let $V$ be a bounded $\vbbK$-fold algebra satisfying skew-symmetry for
$\vT\in\End(V)\even^2$ and $M$ a bounded $\vbbK$-fold $V$-module with
translation endomorphism $\vT$. 
If $a, b\in V$ and $b, a$ are both dual in the exchange channel
then $a, b$ are local in the direct channel.
\end{prop}

The idea of the proof is to argue that $a(bc)=(ab)c=(ba)c=b(ac)$.

\begin{proof}
Let $c\in M$. 
We may assume that equations \eqref{E:dual ab} and
\begin{equation}
\label{E:dual ba}
b(\vx+\vz)a(\vz)c
\; =\;
\sum_k\frac{e^k(\vz,\vx)}{(\vx+\vz)^{\vh_{bc}^k}},
\qquad
(b(\vx)a)(\vz)c
\; =\;
\sum_k\frac{e^k(\vz,\vx)}{(\vz+\vx)^{\vh_{bc}^k}}
\end{equation}
are satisfied where 
\begin{equation}
\label{E:duality coeff ab}
d^j(\vw,\vx)
\; =\;
\sum_{i,k}\:
\vw^{-\vh_{bc}^k}\,\vx^{-\vh_{ab}^i}\:
p_{ijk}^{ab}(\vw,\vx),
\end{equation}
\begin{equation}
\label{E:duality coeff ba}
e^k(\vz,\vx)
\; =\;
\sum_{i,j}\:
\vz^{-\vh_{ac}^j}\,\vx^{-\vh_{ab}^i}\:
p_{ijk}^{ba}(\vz,\vx),
\end{equation}
$p_{ijk}^{ab}(\vz,\vw), p_{ijk}^{ba}(\vz,\vw)\in M\pau{\vz,\vw}$,
and $\vh_p^l\in\vbbK$ such that 
$\vh_p^l\notin\vh_p^{l'}+\Z^2$ for $l\ne l'$ and $p\in\set{ab,ac,bc}$.

Inserting \eqref{E:duality coeff ab} into \eqref{E:dual ab} we get
$$
(a(\vx)b)(\vw)c
\; =\;
\sum_{i,j,k}\:
(\vw+\vx)^{-\vh_{ac}^j}\,\vw^{-\vh_{bc}^k}\,\vx^{-\vh_{ab}^i}\:
p_{ijk}^{ab}(\vw,\vx).
$$
On the other hand,
applying skew-symmetry and 
inserting \eqref{E:duality coeff ba} into \eqref{E:dual ba} we get
\begin{align}
\notag
(a(\vx)b)(\vw)c
\; =\;
&
\paraab\,
(e^{\vx\vT}b(-\vx)a)(\vw)c
\\
\notag
\; =\;
&
\paraab\,
e^{\vx\del_{\vw}}(b(-\vx)a)(\vw)c
\\
\notag
\; =\;
&
\paraab\,
e^{\vx\del_{\vw}}
\sum_{i,j,k}\:
\vw^{-\vh_{ac}^j}\,(\vw-\vx)^{-\vh_{bc}^k}\,(-\vx)^{-\vh_{ab}^i}\:
p_{ijk}^{ba}(\vw,-\vx)
\\
\notag
\; =\;
&
\paraab
\sum_{i,j,k}\:
(\vw+\vx)^{-\vh_{ac}^j}\,\vw^{-\vh_{bc}^k}\,(-\vx)^{-\vh_{ab}^i}\:
p_{ijk}^{ba}(\vw+\vx,-\vx).
\end{align}
Thus Proposition \ref{P:ko}\,i) yields
\begin{equation}
\label{E:pijkab}
p_{ijk}^{ab}(\vw,\vx)
\; =\;
\paraab\,
(-1)^{\vh_{ab}^i}\; p_{ijk}^{ba}(\vw+\vx,-\vx).
\end{equation}

Inserting \eqref{E:duality coeff ab} into \eqref{E:dual ab} we get
\begin{align}
\notag
a(\vz)b(\vw)c
\; =\;
&e^{-\vw\del_{\vz}}
a(\vz+\vw)b(\vw)c
\\
\notag
\; =\;
&
e^{-\vw\del_{\vz}}
\sum_{i,j,k}\:
(\vz+\vw)^{-\vh_{ac}^j}\,\vw^{-\vh_{bc}^k}\,\vz^{-\vh_{ab}^i}\:
p_{ijk}^{ab}(\vw,\vz)
\\
\notag
\; =\;
&
\sum_{i,j,k}\:
\vz^{-\vh_{ac}^j}\,\vw^{-\vh_{bc}^k}\,(\vz-\vw)^{-\vh_{ab}^i}\:
p_{ijk}^{ab}(\vw,\vz-\vw).
\end{align}
In the same way,
inserting \eqref{E:duality coeff ba} into \eqref{E:dual ba} we obtain
$$
b(\vw)a(\vz)c
\; =\;
\sum_{i,j,k}\:
\vz^{-\vh_{ac}^j}\,\vw^{-\vh_{bc}^k}\,(\vw-\vz)^{-\vh_{ab}^i}\:
p_{ijk}^{ba}(\vz,\vw-\vz).
$$
Because $p_{ijk}^{ab}$ and $p_{ijk}^{ba}$ are power series 
\eqref{E:pijkab} implies
$$
p_{ijk}^{ab}(\vw,\vz-\vw)
\; =\;
\paraab\,
(-1)^{\vh_{ab}^i}\; p_{ijk}^{ba}(\vz,\vw-\vz).
$$
Thus locality holds.
\end{proof}

\section{Modules}
\label{S:mod}

We define the notion of a module over an OPE-algebra $V$ and prove that
if $V$ is uniformly local then $V$ is a $V$-module.

\medskip

Let $V$ be a bounded $\vbbK$-fold algebra and $M$ a
bounded $\vbbK$-fold $V$-module.
Then $M$ is {\it dual} if there exist 
$r_{ab}\in\N$ and $\vh_{ab}^i\in\vbbK$ for any
$a\in V, b\in V\cup M, i\in [r_{ab}]$ such that
for any $a, b\in V, c\in M$ there exist
$p_{ijk}(\vw,\vx)\in M\pau{\vw,\vx}$ for $i\in [r_{ab}], j\in [r_{ac}], 
k\in [r_{bc}]$ such that
\begin{equation}
\label{E:dualp}
a(\vx+\vw)b(\vw)c
\; =\;
\sum_{i,j,k}\:
\frac{p_{ijk}(\vw,\vx)}
{\vx^{\vh_{ab}^i}\,(\vx+\vw)^{\vh_{ac}^j}\,\vw^{\vh_{bc}^k}},
\qquad
(a(\vx)b)(\vw)c
\; =\;
\sum_{i,j,k}\:
\frac{p_{ijk}(\vw,\vx)}
{\vx^{\vh_{ab}^i}\,(\vw+\vx)^{\vh_{ac}^j}\,\vw^{\vh_{bc}^k}}.
\end{equation}
If $M$ is dual then duality in the direct and in the exchange channel
are satisfied.

Propositions \ref{P:dual skew} and \ref{P:skew} show that
if $V$ is an OPE-algebra, $\vT$ is a translation endomorphism of $M$, 
and $M$ is dual then $V$ is local on $M$.
Thus the following definition makes sense.

\begin{defi}
Let $V$ be an OPE-algebra.
A $V$-{\it module} is a bounded $\vbbK$-fold $V$-module $M$ such that
$\vT$ is a translation endomorphism of $M$, $M$ is dual,  
and $Y: V\to\cF(M)$ is a $\K^2$-fold algebra morphism.
\end{defi}

\begin{lem}
Let $c\upal\in E\sqbrack{\vz_1, \dots, \vz_r}$ and $\vh_{ij}^k\in\K^2$
for any $i, j\in [r], k\in [r_{ij}]$, and $\al\in\prod_{i<j}[r_{ij}]$
such that $\vh_{ij}^k\notin\vh_{ij}^l+\Z^2$ for $k\ne l$ and 
$$
\sum_{\al}
\frac{c\upal(\vz_1, \dots, \vz_r)}
{\prod_{i<j}(\vz_i-\vz_j)^{\vh_{ij}^{\al(i,j)}}}
\; =\;
0.
$$
Then $c\upal=0$ for any $\al$.
\end{lem}

\begin{proof}
From $(\vz-\vw)^{\vh}\in\K[\vw^{\K}]\set{\vz}$ we get
$$
\prod_{(i,j)\ne (i',j')}(\vz_i-\vz_j)^{\vh_{ij}^{\al(i,j)}}
\;\in\;
\K\sqbrack{\vz_1}\dots\sqbrack{\vz_{i'},\vz_{j'}}\set{\vz_{i'+1}}\dots
\widehat{\set{\vz_{j'}}}\dots\set{\vz_r}
$$
for any $i'<j'$. 
Thus the claim follows from Proposition \ref{P:ko}\,i)
by induction on the number of pairs $i<j$.
\end{proof}

\begin{prop}
\label{P:sing}
Let $V$ be a uniformly local OPE-algebra, $a, a^i\in V$, and 
\eqref{E:mult loc} an OPE of $a^1(\vz), \dots, \linebreak[0] a^r(\vz)$.
Then
$$
c\upal(\vz_1, \dots, \vz_r)a
\; \in\;
\sum_{\ka}
\vz_1^{-\vh_{a^1 a}^{\ka 1}}\dots \vz_r^{-\vh_{a^r a}^{\ka r}}
V\pau{\vz_1, \dots, \vz_r}
$$
for any $\al$ where $\ka\in\prod_{i=1}^r [r_{a^i a}]$.
\end{prop}

\begin{proof}
We may assume that the OPE \eqref{E:mult loc} is reduced.
Consider multiple locality for $a^1(\vz_1), \dots, a^r(\vz_r), \linebreak[0]
a(\vz)$ and $\si=1$.
Apply it to $1\in V$ and set $\vz=0$. 
By assumption, 
$c^{(\al,\ka)}(\vz_1, \dots, \vz_r, 0)1\in V\pau{\vz_1, \dots, \vz_r}$.

On the other hand, 
consider multiple locality for $a^1(\vz_1), \dots, a^r(\vz_r)$ and $\si=1$.
Apply it to $a$. In both cases we get $a^1(\vz_1)\dots a^r(\vz_r)a$
on the left-hand side. Equating the right-hand sides and applying the Lemma
yields the claim.
\end{proof}

\begin{cor}
Let $V$ be a uniformly local OPE-algebra.
Then $V$ is a $V$-module.
\end{cor}

\begin{proof}
The fact that $V$ is dual follows from Proposition \ref{P:sing} and from
Proposition \ref{P:dual loc} and its proof.
That $Y$ is a morphism follows from the theorem in section \ref{S:mult loc}.
\end{proof}

\begin{acknowledgements}
This paper consists mostly of results of my thesis \cite{rosellen.phd}
written under the supervision of Yu.~Manin and W.~Nahm.
My work on OPE-algebras
originated from a collaboration with D.~Huybrechts and M.~Lehn.
D.~Orlov showed me how to prove that $c^i(\vz,\vw)1$ are power series.
D.~van Straten and A.~Matsuo invited me to give talks about OPE-algebras.
My post-doc at Stockholm is made possible by S.~Merkulov.
I am grateful to all of them.
\end{acknowledgements}

\bibliographystyle{amsalpha}
\bibliography{lit/bib,lit/collections}

\end{document}